\documentclass[openany]{book}

\usepackage[medium]{dgruyter}

\theoremstyle{dgdef}

\theoremstyle{dgthm}
\newtheorem{theorem}{Theorem}

\newtheorem{remark}{Remark}

\def\Div{\mbox{div}\,}

\def\bE{{\bf E}}

\def\dfrac#1#2{\displaystyle{#1\over #2}}

\def\bV{{\bf V}}

\begin{document}
\contribution


  \contributionauthor*[1]{Olga Rozanova}			

  \runningauthor{O. Rozanova}
  \affil[1]{Mathematics and Mechanics Department, Lomonosov Moscow State University, Leninskie Gory,
Moscow, 119991, Russian Federation; email: rozanova@mech.math.msu.su}

  \contributiontitle{Regularizing factors for the Euler-Poisson equations}
  \runningtitle{Regularizing factors}
  \abstract{The Cauchy problem for the Euler-Poisson equations without pressure is considered and the question of what additional terms added to the system can delay or completely prevent the loss of smoothness of the solution in a finite time is studied. We review already published and recent results in this field.}

\makecontributiontitle

\section{Introduction}
The  Euler-Poisson equations have the following form
\begin{eqnarray}\label{EP1}
\dfrac{\partial n }{\partial t}& + &\Div(n \bV)=0,\\
\dfrac{\partial \bV }{\partial t} &+& \left( \bV \cdot \nabla \right)
\bV + {\alpha \frac{\nabla p(n)}{n}} =\,k \,  \nabla \Phi - {\nu \bV}  +{\mu \Delta \bV}, \label{EP2}\\ \Delta \Phi &=&n-n_0,\label{EP3}
\end{eqnarray}
where $n$ is the density, $\Phi$ is a force potential,  $\bV$ is the vector of velocity, $n_0={\rm const}\ge 0$ is the density background,
$k=\rm const$. Components of the solution depend on $t>0$ and $x\in {\mathbb R}^n$.

If the equations describe a medium consisting of electrons (the case of plasma or semi-conductors), then the force is repulsive, i.e. $k>0$, if the medium consists of gravitating particles as in the astrophysics models, then the force is attractive, i.e.
$k<0$. Further, $\alpha  p(n)\ge 0$ is the pressure, $- {\nu \bV}$ and   ${\mu \Delta \bV}$ are the friction and viscosity terms.
Here $\alpha, \nu, \mu $ are nonnegative coefficients. If we set the coefficient to zero, then we do not consider the corresponding term in the model. We call the system with $\alpha =\nu =\mu= 0 $ {\it the original Euler-Poisson equations}.

It is well known that solutions to the Cauchy problem for the original Euler-Poisson equations are prone to loss of smoothness (e.g. \cite{Rozanova:ELT}). Nevertheless, if we use the Euler-Poisson equations for modeling physical processes, the nature of processes sometimes dictates that the equations are valid only for smooth solutions. For example, the cold (electron) plasma can be reasonably modeled by the original Euler-Poisson equations in the repulsive case, but only until the moment of a singularity formation. Singularities are associated with the formation of the delta function in the component of density, and this process leads to the heating of the plasma. Therefore we have to move to another model (adding some terms).

Thus, our main problem is to study the influence of additional terms on the lifetime of a smooth solution to the Cauchy problem for the Euler-Poisson system. Since for the multidimensional case, the results about the possibility of the existence of globally smooth solutions are very scarce, we restrict ourselves to the one-dimensional case. In what follows we focus on the repulsive case with a nonzero density background corresponding to a cold plasma. This is the most interesting and complex case due to the oscillatory nature of the solution and the most important model from a practical point of view.

For the original Euler-Poisson system, there is a criterion for the formation of a singularity from the initial data, that is, the exact class of initial data that gives a globally smooth solution is known (we call it the class of smoothness).

Let us present a short list of known results for the case of ''pure`` factors, when only one of the coefficients $\alpha, \nu, \mu$ is nonzero.

\begin{itemize}
\item $\nu>0$,  $\alpha =\mu= 0 $ (the friction, Sec.3)
\begin{itemize}
\item $\nu={\rm const >0}$: the class of smoothness enlarges with $\nu$, but for any $\nu$ there exist initial data such that the respective solution blows up in a finite time;

    \item $\nu=\nu (n)$: there is a dependence on $n$ such that the solution preserves global smoothness for all initial data;
 \end{itemize}
 \item $\alpha>0$,  $\mu =\nu= 0 $ (the pressure, Sec.4):  the class of smoothness does not enlarge, but the type of singularity changes (weakens);
\item $\mu>0$,  $\alpha =\nu= 0 $ (the viscosity, Sec.5)
\begin{itemize}
\item $\mu={\rm const >0}$:  the solution keeps the global smoothness for all initial data;
   \item an exotic viscous term depending on $n$ and $\bV$ may do not prevent  blow-up.
   \end{itemize}
 \end{itemize}
We do not consider here the combination of factors (for example, friction and pressure), but we believe that in the previous context, it is not a difficult problem.

 We can also consider a stochastic counterpart of the original Euler-Poisson equations and then study the smoothness of deterministic characteristics of the stochastic process (the probability density and analogs of velocity and force acting between particles), Sec.6. 
 These deterministic characteristics are globally smooth in $t$ and tend to the solution of the original Euler-Poisson system as the parameter of the stochastic perturbation tends to zero.


\section{Euler-Poisson equations in alternative form}

Let us transform \eqref{EP1} -- \eqref{EP3} to a more convenient form. To do this, we introduce the vector function
$\bE=-  \nabla \Phi$. From \eqref{EP3} we get
\begin{eqnarray}
n=n_0- \Div \bE, \label{n}\end{eqnarray}
which allows us to eliminate $ n $ by substituting \eqref{n} into \eqref{EP1} and obtain $\Div {\bf X}=0$, where ${\bf X}=\frac{\partial \bE }{ \partial t} + \bV \Div \bE
- n_0 \bV$. According to Helmholtz's theorem, ${\bf X}=\nabla F + {\rm rot}\,{\bf A} $, where $F$ and $\bf A$ are scalar and vector potentials, respectively. Since ${\rm rot}\,{\bE}=0$, then in the case when ${\rm rot} (n \bV)=0 $ (in the one-dimensional case in space this is obviously true)
we have ${\rm rot}\,{\bf X}=0$. Assuming that the components of the vector $\bf X$ decrease quickly enough as $|x|\to \infty$, according to the theorem on restoring a smooth vector field from its divergence and rotor \cite{Rozanova:Kochin} we obtain that $\bf X $ is unique and equal to zero.

The resulting system  is
\begin{eqnarray}\label{4}
\dfrac{\partial \bV }{\partial t} + \left( \bV \cdot \nabla \right)
\bV &=& \, - k \bE {{ - \nu \bV + \mu \Delta \bV - \alpha  \frac{\nabla p(n)}{n}}},\\ \frac{\partial \bE }{\partial t} +  \bV \Div \bE
& =& n_0\bV.
\end{eqnarray}

In what follows we consider the 1D case, $k=1$, $n_0=1$ and denote $\bV=V$, $\bE=E$.

\subsection{Original system, $\nu=\mu=\alpha =0$ }
For this case, we obtain a non-strictly hyperbolic system
\begin{eqnarray}\label{q11}
\dfrac{\partial V }{\partial t}  + V \dfrac{\partial V}{\partial x}= - E,\qquad
\dfrac{\partial E }{\partial t} +  V \dfrac{\partial E}{\partial x} = V,
\end{eqnarray}
which we consider together with initial data
\begin{eqnarray}
(V,E)|_{t=0}=(V_0,E_0).\label{q1CD}
\end{eqnarray}
Let us denote $(V_x,E_x)=(v,e)$. Then
along characteristics starting point $x_0\in \mathbb R$ we have the following dynamics:
\begin{eqnarray}\label{q1}
\dot V = -E, \qquad \dot E = V,\qquad
  \dot v= -e-v^2, \qquad \dot e = v (1 - e),
 \end{eqnarray}
\begin{equation*}\label{eq_cd}
  (V,E,v,e)|_{x=x_0, t=0}=(V_0(x_0), E_0(x_0), v_0(x_0), e_0(x_0)).
 \end{equation*}

\begin{theorem}\cite{Rozanova:RChZAMP}
For the existence and uniqueness of continuously differentiable  $ 2 \pi - $ periodic in time  solution $ (V,E)$ of \eqref{q11}, \eqref{q1CD} (where \eqref{q1CD} belong to $C^2({\mathbb R})$ ) it is necessary and sufficient  that at each point $ x \in  \mathbb R$ 
\begin{eqnarray} \label{crit2}
\Delta(x)=v_0^2(x)+2e_0(x)-1<0.
\end{eqnarray}
\end{theorem}

Fig.\ref{Pic4} presents the domain of smoothness on the plane $(v_0, e_0)=(v_0(x_0), e_0(x_0))$, such the solution of \eqref{q1}, \eqref{q1CD} does not blow up along the  characteristics starting from $x_0$  (left, bounded by a dotted line).

\section{Influence of friction ( $\alpha =\mu=0$)}

\subsection{Constant coefficient of friction }
\begin{center}
\begin{figure}[h]
\begin{minipage}{0.45\columnwidth}
\centerline{
\includegraphics[width=0.9\columnwidth]{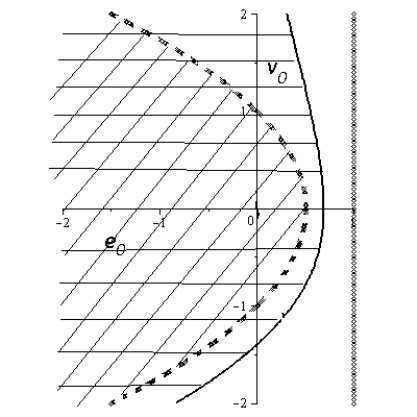}}
\end{minipage}
\hspace{1cm}
\begin{minipage}{0.45\columnwidth}
\centerline{
\includegraphics[width=0.9\columnwidth]{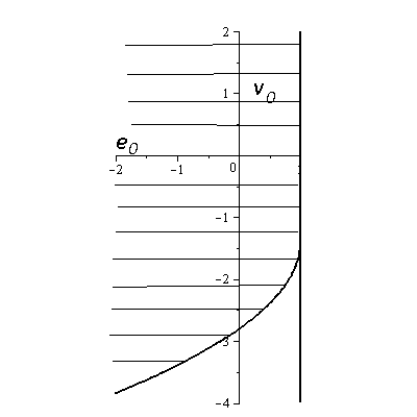}}
\end{minipage}
\caption{Left: $ \nu < 2 $: the domain of smoothness for  $ 0< \nu < 2 $ (shaded) in comparison with the domain of smoothness for $ \nu = 0 $ (double shaded). Right: the domain of smoothness for  $ \nu >2$ (shaded).}\label{Pic4}
\end{figure}
\end{center}

The case  $\nu= {\rm const} >0$ was considered in \cite{Rozanova:ELT} in terms on the Euler-Poisson equations and in \cite{Rozanova:RChD_fric} in terms of non strictly hyperbolic system
\begin{eqnarray*}
\dfrac{\partial V }{\partial t}  + V \dfrac{\partial V}{\partial x}= - E-\nu V,\qquad
\dfrac{\partial E }{\partial t}+  V \dfrac{\partial E}{\partial x} = V.
\end{eqnarray*}
 The characteristic system  is
\begin{eqnarray}\label{q2}
\dot V = -E-\nu V, \qquad \dot E = V, \qquad \dot v= -e-v^2-
\nu v, \qquad \dot e = v (1 - e).
\end{eqnarray}
The behavior of the solution depends on the intensity of friction $\nu$. Namely, the analysis of the phase plane $(e,v)$ for the characteristic system \eqref{q2} shows that
\begin{itemize}
 \item  if $\nu=0$, then we have one equilibrium point $ (0,0) $, a center (the original Euler-Poisson system), the solution is oscillatory and periodic in $t$;
 \item if $0<\nu<2$, then  we have one equilibrium point $(0,0)$, a stable focus, if the solution does not blow up, it is oscillatory and the amplitude decays to zero;
 \item  if $\nu>2$, then we have  three equilibria, $(0,0)$, a  stable node,  $(1,-\frac12(\nu-\sqrt{\nu^2-4})$, a saddle,  $(1,-\frac12(\nu+\sqrt{\nu^2-4})$, an  unstable node (if $\nu=0$, then the saddle and the unstable node merge into a saddle-node),
     if the solution does not blow up and decays to zero without oscillations;
      \end{itemize}
      It is possible to find the explicit expression for the curve, separating the class of smoothness from the class of blow-up on the phase plane $(v_0, e_0)$ (see \cite{Rozanova:RChD_fric}), however we do not write here this quite cumbersome formula (the curve is composed of three phase curves). For greater clarity, we present Pic.\ref{Pic4}, which shows the extension of the class of smoothness with $\nu$ and the change of character of this curve for large $\nu$. In particular, for any fixed initial data we can find $\nu={\rm const}>0$ such that the solution to the Cauchy problem remains smooth for all $t>0$. However, for any fixed arbitrary large $\nu={\rm const}>0$ one can find initial data such that the solution to the Cauchy problem blows up within a finite time.
      
\subsection{Density dependent friction: $\gamma=\mu=0$,  $\nu = \nu(n)>0$}
Now we assume a dependence of the friction coefficient on $n$, which is more realistic for the models of a cold plasma. Let us consider the following problem with the initial data from the class of analytical functions, denoted below as ${\mathcal A}$ (a technical requirement):
\begin{equation*}\label{2}
\dfrac{\partial V}{\partial t}+ V  \dfrac{\partial V}{\partial x}=-E
-\nu(n) V, \quad \dfrac{\partial E}{\partial t}+ V  \dfrac{\partial
E}{\partial x}=V, \quad n=1- \dfrac{\partial E}{\partial x},
\end{equation*}
\begin{equation*}\label{CD}
(V,E)|_{t=0}=(V_0(x), E_0(x))\in {\mathcal A} ({\mathbb R}).
\end{equation*}

One can prove that there exists a choice of $\nu(n)$ such that the solution remains smooth for all possible choices of smooth initial data.

\begin{theorem} \cite{Rozanova:R_PD} Let
$f(n)\in {\mathcal A}({\mathbb R}_+) $ be a non-negative function
satisfying conditions
\begin{equation*}\label{aid_cond2}
    \lim\limits_{\eta\to\infty}\dfrac{\eta f'(\eta)}{f(\eta)}={\rm
    const}, \qquad\qquad
\int\limits_{\eta_0>0}^{+\infty}\dfrac{ f(\eta)}{\eta^2}\,
d\eta=\infty,
\end{equation*}
$\nu(n) = \epsilon f (n)$, $\epsilon =\rm
const >0$ is a small parameter.
Under the assumption that the formation of singularity is associated with a gradient catastrophe (unboundedness of the first derivatives), the problem admits a global in time classical ($C^1$-smooth) solution. Otherwise, one can find the
data such that derivatives of the solution blow up in a finite time.
\end{theorem}

The prototypic function is
$\nu(n)=\nu_0 n^\gamma$ with  the threshold value $\gamma=1$. In other words,
for
$\gamma\ge 1$ the solution does not form
a gradient catastrophe.

\section{Influence of the pressure, $\nu=\mu = 0$, $\alpha \ne 0$}

The study of the influence of the pressure of the Euler-Poisson system is highly nontrivial since the system became strictly hyperbolic
in terms of $n$ and $V$.  Here we only announce the final result and skip quite cumbersome proof that is now in the preparation for publication.

\begin{theorem} \label{T1}  Assume $p(n)=\frac{1}{\gamma} n^\gamma$, $\gamma>1$. A continuously differentiable  $ 2 \pi- $ periodic in time  solution $ (V,E)$ of
\begin{eqnarray*}\label{q1p}
&&\dfrac{\partial V }{\partial t}  + V \dfrac{\partial V}{\partial x}= - E - \alpha \frac{1}{n} \dfrac{\partial p(n)}{\partial x},\quad
\dfrac{\partial E }{\partial t} +  V \dfrac{\partial E}{\partial x} = V, \quad
n=1-\dfrac{\partial E }{\partial x},\\
&&(V,E)|_{t=0}=(V_0,E_0)\in C^2({\mathbb R}), \quad n_0=1-E'_0.
\end{eqnarray*}
 exists if and only if
\begin{equation} \label{crit2p}
\Delta_p (x)=v_0^2 (x)+2e_0(x)-1 <\alpha \frac{ e_0'(x)^2}{(1-e_0(x))^{3-\gamma}}
\end{equation}
holds at each point $ x \in  \mathbb R$.
\end{theorem}

We can see that the pressure generally does not remove or postpone a singularity. Indeed, let us consider the initial data in the form of a standard laser pulse with $\Phi_0= a \exp(-x^2), \, a>0$, i.e.
$$
V_0=0, \quad E_0= a (\exp(-x^2))'.
$$
For $\alpha=0$ the most dangerous point for the blow-up is $x=0$, since $\Delta$ takes its maximum value there. But at this moment
$E_0''(x)=0$, so $\Delta_p(0)=\Delta(0)$ and the requirements for the intensity $a$ to satisfy the conditions \eqref{crit2p} and \eqref{crit2} are the same.
 However, it is important that the type of singularity changes. Namely, if for $\alpha=0$ the formation of singularity implies the gradient catastrophe for $V$ and $E$, and the strong singularity for $n$,
for $\alpha>0$   the formation of singularity implies the gradient catastrophe for $V$ and $n$, and the weak singularity for $E$.

\subsection{Example: $\gamma=2$}
In this case the system \eqref{q1p} takes the form
\begin{eqnarray*}
\dfrac{\partial V }{\partial t}  + V \dfrac{\partial V}{\partial x}= - E + \alpha \dfrac{\partial^2 E}{\partial x^2},\quad
\dfrac{\partial E }{\partial t} +  V \dfrac{\partial E}{\partial x} = V,
\end{eqnarray*}
It is a particular case of
\begin{equation}\label{msp}
\dfrac{\partial \mathfrak V }{\partial t} + V_1\dfrac{\partial\mathfrak V }{\partial x}=  Q \mathfrak V + B \dfrac{\partial^2\mathfrak V }{\partial x^2},
\end{equation}
see \cite{Rozanova:R_nonstr}, where
$\mathfrak V=(V_1, V_2, \dots, V_n)$, $V_i=V_i(t,x)$,  $Q$ and $B$  are $n\times n$ constant matrices. Here $V_1=V$, $V_2=E$, $Q= \left(\begin{array}{cc}0 & -1\\ 1& 0 \end{array}\right)$, $B= \left(\begin{array}{cc}0 & \alpha\\ 0& 0 \end{array}\right)$.
Let us notice that $B \dfrac{\partial^2\mathfrak V }{\partial x^2}$ looks like a viscous term, however, in fact, its sense is different.

\section{Influence of viscosity, $\nu=\alpha = 0$, $\mu={\rm const}> 0$}

Now we consider the following system:
\begin{equation}\label{q1mu}
\dfrac{\partial V }{\partial t}  + V \dfrac{\partial V}{\partial x}=- E +\mu \dfrac{\partial^2 V}{\partial x^2},\qquad
\dfrac{\partial E }{\partial t} + V \dfrac{\partial E}{\partial x} = V,
\end{equation}
a particular case of \eqref{msp} with  $B= \left(\begin{array}{cc}\mu & 0\\ 0& 0 \end{array}\right)$.

\begin{theorem} \label{T4}  Assume $\mu={\rm const}> 0$. Then the solution to the Cauchy problem \eqref{q1mu}, \eqref{q1CD}, with the data from the class $C^1({\mathbb R})$ keeps this smoothness for all $t>0$.
 \end{theorem}

{\it Sketch of the proof.} We apply the Cole-Hopf transformation, e.g. \cite{Rozanova:CH}, to the first equation of \eqref{q1mu}: $U_x=V$, $U=- 2 \mu \ln z,$ and obtain 
\begin{equation*}
z_t   =\mu  z_{xx}+ \Psi z, \qquad \Psi_x=-E.
\end{equation*}
If $\Psi\in {\mathbb H}^{l,\frac{l}{2}}({\mathbb R}\times (0,T))$, $z(0,x)\in {\mathbb H}^{l+2}({\mathbb R})$, then $z\in {\mathbb H}^{l+2,\frac{l}{2}+1}(({\mathbb R}\times (0,T))$, $l\ge 0$, where ${\mathbb H}^{l,\frac{l}{2}}$ is the notation for the parabolic H\"older space  \cite{Rozanova:LSU}. Assume that for some time $t_*>0$ the component $E(t,x)$ forms a jump in $x$. However, the smoothness properties of $\Psi$ are better, $\Psi\in{C}^{0,0}(({\mathbb R}\times (0,T)) $ and therefore $z\in {C}^{2,1}(({\mathbb R}\times (0,T))$ and
 $V\in {C}^{1,1}(({\mathbb R}\times (0,T))$. From the second equation of \eqref{q1mu}, hyperbolic and linear with respect to $E$,
 we see that $E$ has the same smoothness as $V$, so we obtain a contradiction.

\begin{remark} Let us notice that the viscosity does not always prevent the blow-up. For example, in \cite{Rozanova:ELT} the system
\begin{equation*}\label{q1}
V_t  + V V_x=  E +\mu \left(\frac{ V_x}{n}\right)_x,\qquad E_t + V E_x = n_0 V, \qquad
 E_x = n_0-n.
\end{equation*}
was considered and a condition in terms of the initial data for the finite time blow-up was found, such that
$$\frac{V_x}{n}\to -\infty, \qquad  t\to T_*<\infty.  $$
\end{remark}
\begin{remark} In the model of perturbation of a rest state in the Reyleigh-B\'enard convection  theory \cite{Rozanova:Drazin2002} and
the model of a stratified fluid near a rest state  in a gravitational field \cite{Rozanova:Baidulov} there arises the system
\eqref{msp} with  $B= \left(\begin{array}{cc}\mu & 0\\ 0& \kappa \end{array}\right)$, $\mu, \kappa={\rm const} >0$. Thus, the viscosity matrix here is non-degenerate, and one can expect the $C^\infty$ -- smoothness of the solution after this parabolic regularization. Indeed, we can apply the same method as in the proof of Theorem \ref{T4}. In this new case, the analog of the second equation in \eqref{q1mu} is parabolic, and the smoothness of $E$ increases compared to $V$. This results in a greater smoothness of $\Phi$, and we can introduce an iteration process that allows us to prove that the solution belongs to $C^k$ for any $k\in \mathbb{N}$.
The global in $t$ existence of the classical solution to the Cauchy problem follows from results of \cite{Rozanova:LSU}.
\end{remark}

\section{Stochastic regularization}

Let us give the idea of a probabilistic method, which can be applied to any non-strictly hyperbolic analog of  \eqref{msp},
\begin{equation}
\label{base11}
\begin{array}{l}
\dfrac{\partial \bV }{\partial t} +  V_1 \mathbb E \dfrac{\partial\bV }{\partial x}=  Q \bV.
\end{array}
\end{equation}
The method consists of a stochastic perturbation along characteristics \cite{Rozanova:AKR}.

Namely, let us consider a stochastic counterpart of \eqref{base11}
\begin{eqnarray*}
d X(t)&=& \mathfrak V_1(t) dt + \sigma dW,\\
d \mathfrak V (t)& =& Q \mathfrak V (t),
\end{eqnarray*}
where $W$ is a standard Wiener process, $\sigma>0$, $\mathfrak V (t)=\bV (t, X(t))$, $\mathfrak V_1 (t)=V_1 (t, X(t))$.

Let $P(t,x, v_1,..,v_n)$ be joint probability function for processes $(X(t),\mathfrak V (t))$, it satisfies to the Fokker-Planck equation (it is standard and we do not write it here) with the initial data
$$P_0(t,x, v_1,..,v_n)=\prod\limits_{i=1}^n \delta(v_i-V_i(0,x))\, f_0(x), \quad f_0\in L_1(\mathbb R).$$

Let  us introduce
\begin{eqnarray*}
\rho(t,x) &=&
\int\limits_{\mathbb R^n} P(t,x,v_1,...,v_n)\, dv_1...dv_n
\\
\hat\bV(t,x)&=& \frac {1}{\rho(t,x)}\int\limits_{\mathbb R^n} v P(t,x,v_1,...,v_n)\, dv_1...dv_n,\quad v=(v_1,...,v_n).
\end{eqnarray*}

Then  we obtain the system
 \begin{eqnarray*}
 \dfrac{\partial \rho  }{\partial t} &+&   \dfrac{\partial \rho \hat V_1 }{\partial x}=  \frac{\sigma^2}{2}
\dfrac{\partial^2 \rho  }{\partial x^2},\\
 \dfrac{\partial \rho \hat\bV }{\partial t} &+&   \dfrac{\partial \rho \hat V_1 \hat\bV }{\partial x}=  Q \rho \hat\bV+ \frac{\sigma^2}{2}
\dfrac{\partial^2 \rho  \hat\bV }{\partial x^2}-\int\limits_{\mathbb R} \, (v-\hat \bV) (v_1-\hat V_1) P_x \, d v,
\end{eqnarray*}
One can prove that
$
\hat\bV \to \bV,$ $\sigma\to 0,$ $t>0, $ $ x\in\mathbb R,$
for a {\it continuous} $\bV$.

\subsection{Example: regularization of the original equations of cold plasma}

We apply the method to \eqref{q11}, it is a particular case of \eqref{base11}. 
The result of the stochastic regularization is
 \begin{eqnarray*}
 \dfrac{\partial \rho  }{\partial t} &+&   \dfrac{\partial \rho \hat V }{\partial x}=  \frac{\sigma^2}{2}
\dfrac{\partial^2 \rho  }{\partial x^2},\\
 \dfrac{\partial \rho \hat V }{\partial t} &+&   \dfrac{\partial \rho \hat V^2  }{\partial x}=  - \rho \hat E+ \frac{\sigma^2}{2}
\dfrac{\partial^2 \rho  \hat V }{\partial x^2}-\int\limits_{\mathbb R^2} \, (v-\hat V)^2  P_x \, d v d e,\\
 \dfrac{\partial \rho \hat E }{\partial t} &+&   \dfrac{\partial \rho \hat V  \hat E }{\partial x}=  \rho \hat V+ \frac{\sigma^2}{2}
\dfrac{\partial^2 \rho  \hat E }{\partial x^2}-\int\limits_{\mathbb R^2} \, (v-\hat V) (e-\hat E) P_x \, dv de,
\end{eqnarray*}
Here the components of solution is an artificial density $\rho$ and $\hat V$, $\hat E$, regularized analogs of $V, E$, namely,
\begin{eqnarray*}
&\rho(t,x) =
\int\limits_{\mathbb R^2} P(t,x,v,e)\, dv de,
\\
&\hat V(t,x)= \frac {1}{\rho(t,x)}\int\limits_{\mathbb R^2} v P(t,x,v,e)\, dv de,\qquad
\hat E(t,x)= \frac {1}{\rho(t,x)}\int\limits_{\mathbb R^2} e P(t,x,v,e)\, dv de.
\end{eqnarray*}

Thus, based on the original system \eqref{q11}, solutions of which can blow under certain conditions, we obtain a new system, with globally in-time smooth solutions. For small $\sigma>0$ the solution of this new system is similar to the solution of the original system until the blow-up moment, however after this moment the behavior of the smooth solution is different. In particular, the delta-shaped waves do not form in the component of density. We do not present the proof here, which is very similar to \cite{Rozanova:AKR}.

\begin{funding}
  This work was financially supported by the Russian Science Foundation under grant no. 23-11-00056.
\end{funding}

\end{document}